\documentclass{amsart}
\usepackage{amsmath,amsfonts,amsthm}
\usepackage{amssymb,latexsym}
\usepackage{graphics}
\usepackage{rotating}
\usepackage{float}
\usepackage[colorlinks]{hyperref}

 \usepackage[all]{xypic}
\usepackage{amssymb}

\theoremstyle{plain}
\theoremstyle{definition}
\newtheorem{theorem}{Theorem}
\newtheorem{lemma}[theorem]{Lemma}

\newtheorem{example}[theorem]{Example}

\theoremstyle{remark}

\numberwithin{equation}{section}

\title{Residual fertility and delay in sterile insect population dynamics}
\author[Ruqaya Hussein, Sergey Kryzhevich, Khosro Tajbakhsh]{Ruqaya Hussein, Sergey Kryzhevich, Khosro Tajbakhsh}
\begin{document}
\begin{abstract}
The sterile insect technique controls mosquito-borne diseases such as malaria, dengue, and yellow fever through either eradication or depressing the associated vector population. We formulate a three-dimensional delayed mosquito population suppression model with a saturated release rate to explore the interactive dynamics between wild, sterile, and non-sterile mosquitoes, focusing on the delay and residual fertility in the interactive dynamics among insects. We investigate the stability of the positive equilibrium and derive the Hopf bifurcation conditions. We establish the stability conditions for the positive equilibrium and examine how the time delay ($\tau$) and residual fertility affect the non-sterile insects' dynamics. Below the critical values of the delay, the system remains stable, while beyond that, the Hopf bifurcation is guaranteed under certain circumstances. However, analysis shows a clear band of non-sterile insect population values as residual fertility varies within a very narrow range. This suggests that within this interval, the system exhibits sensitive dependence on the fertility parameter, likely due to underlying nonlinear dynamics. Numerical simulations are presented to support our analytical results, followed by a brief discussion of the findings.
\end{abstract}
\maketitle
\section{Introduction}
\noindent Several methods exist to control mosquito-borne diseases. The traditional approach uses insecticides, but another effective strategy is reducing the wild mosquito population by means of Sterile Insect Technique (SIT), also known as the Sterile Insect Release Method (SIRM). First introduced by E. F. Knipling in 60th 1960s, SIT has been successfully applied in many regions across the world against various insect pests and vectors both in agriculture and medicine. Mosquitoes have four life stages, including egg, larva, pupa, and adult, that take 7 to 20 days from oviposition to adult death on average \cite{bruce}. Their feeding habits are unique as only adult females consume blood (required for egg production), while males feed exclusively on plant nectar. Females target mammals, including humans, as well as birds, reptiles, and amphibians \cite{hard1993}. The principle of SIT against these important insects involves rearing sterile male mosquitoes and releasing them into wild populations. When a wild female mates with a sterile male, she either fails to reproduce or her eggs do not hatch. Over time, this reduces the mosquito population through generations. However, before implementing SIT, optimal release strategies and sterile mosquito quantities must be carefully determined \cite{wang2019}.

Numerous mathematical models have been developed to investigate the interactive dynamics between wild and sterile mosquito populations. Key contributions include studies by but not limited to \cite{cai3,yu,yu2019}. The survey by \cite{cai3} relies on ordinary differential equations (ODEs) and explores three distinct strategies for releasing sterile mosquitoes, including a constant release rate, a release rate proportional to the wild mosquito population, and a saturated release rate proportional to the wild population. Other workers have investigated spreading dynamics by integrating stage structure \cite{ai,lili}, time delays \cite{zh,huan, yu}, spatial diffusion \cite{huan,hua}, and environmental heterogeneity through stochastic equations \cite{huang}.

 Delay differential equations (DDEs) have recently gained increasing attention in mathematical modeling. Introducing delays allows for a more realistic description of biological phenomena and improves predictive accuracy. Such models are also called Retarded Functional Differential Equations (RFDEs). Unlike ODEs, RFDEs are infinite-dimensional dynamical systems, leading to richer dynamics, including oscillations, instability, chaos, and altered stability properties. DDEs are widely used in science, engineering, biology, and medicine for modeling population dynamics, epidemiology, immunology, and neural networks \cite{dubey2000,haque2010}. \cite{cai2} developed a DDE framework to study wild-sterile mosquito interactions and to assess control strategy efficacy. \cite{yu2019} constructed a DDE model accounting for only fertile sterile mosquitoes, focusing on constant-release scenarios, and in a more recent study, \cite{cai2020} evaluated sterile mosquito suppression strategies by incorporating delayed releases. Also, \cite{hui} formulated a delayed suppression model with constant sterile mosquito releases.  

One critical issue in SIT programs is residual fertility in sterile insects that can be challenging when the number of non-sterile sperms/insects in releasing insects/packages exceeds a threshold. This factor is normaly considered as a fixed parameter in SIT modelling; however, in the current work, like \cite{yu} who considred number of released sterile mosquitoes as an independent dynamical equation rather than as a fixed parameter, we incorporate residual fertility as a coupled ODE to wild and sterile insects system of ODEs to analyze consequent dynamics of insects who are escaped from being sterilized enough through irradiation process. Simultaneously, the same delay in releasing both sterile and non-sterile insects as a critical parameter in the success of both SIT and population dynamics of this small portion of sterile insects is studied analytically and numerically. This paper is organized as follows. Section 2 presents the complete formulation of the delay differential equation model. In Section 3, we analyze the stability of the origin, boundary, and positive equilibria, focusing on the circumstances under which Hopf bifurcation can occur in the dynamics, considering the saturated sterile mosquito release rate. Section 4 provides numerical simulations that validate the theoretical results, particularly emphasizing how time delay and residual fertility affect non-sterile insects' population dynamics. Finally, Section 5 concludes by discussing the key findings and their implications. All numerical simulations were performed using MATLAB R2014b, while analytical computations were conducted in Mathematica 12.0. 

\section{The model formulation}
We present a three-dimensional delayed mathematical model with a nonlinear saturated releasing rate containing an ODE governing the dynamics of non-sterile insects coupled to wild and sterile insects' ODEs. Without interactions, the dynamics of the wild mosquito population follow logistic growth. After releasing sterile mosquitoes into the wild population, the system is governed by the following delay differential equations:

\begin{equation}\label{2.1}
\begin{split}
\frac{dw}{dt}&=\Big(\frac{aw(t)+rs(t-\tau)}{1+w(t)+g(t)+s(t-\tau)}-(\xi_{1}(w(t)+g(t)+s(t-\tau))\Big) w(t)\\ 
\frac{dg}{dt}&=\frac{bw(t)}{1+w(t)}-\Big(\xi_{2}(w(t)+g(t)+s(t))\Big) g(t)\\
\frac{ds}{dt}&=\frac{cw(t)}{1+w(t)}+\Big(\frac{rs(t)+aw(t-\tau)}{1+w(t-\tau)+g(t)+s(t)}-(\xi_{3}(w(t-\tau)+g(t)+s(t))\Big) s(t)
\end{split}
\end{equation}
Where
\begin{itemize}
\item  $w(t), g(t)$, and $s(t)$ are populations of wild, sterile, and non-sterile mosquitoes, respectively, at time $t$; 
\item  $N(t)= w(t)+g(t)+s(t)$ is total mosquito population;
\item $a$ is the wild offspring production rate per mating;
\item $\xi_{i}, i=1,2,3$ are density-dependent death rates for wild, sterile, and non-sterile mosquitoes, respectively;
\item $b$ is the sterile mosquito release rate;
\item  $c$ is residual fertility or non-sterile mosquito release rate (small portion of $g$ that remains fertile);
\item $r$ is the non-sterile offspring production rate per mating.
\end{itemize}
The parameters $b$ (sterile release rate) and $c$ (residual fertility proportion) critically influence population dynamics. Adjusting $b$ can suppress wild mosquito numbers, while $r$ and $c$ determine the effectiveness of the SIT program. High $c$ or $r$ may counteract SIT by increasing the reproduction of fertile insects. The coordinates of the system are changed as $W=w-w^*$, $G=g-g^*$, and $S=s-s^*$ to see the positive equilibrium at the origin, and then we expand the Taylor series of the system \eqref{2.2} up to the fourth order. The new equations in terms of new variables W, G, S are as follows.

\begin{equation} \label{2.2}
\begin{split}
\frac{dW}{dt}&= a_{1}+a_{2}W(t)+a_{3}G(t)+a_{4}S(t-\tau)+a_{5} W^2(t)+a_{6} W^3(t)+a_{7} W^4(t)+a_{8} G^2(t)\\&
+a_{9} G^3(t)+a_{10} G^4(t)+a_{11}S^2(t-\tau)+ a_{12}S^3(t-\tau)+h.o.t.,\\
\frac{dG}{dt}&=b_1 W^4(t)+b_2 W^3(t)+b_3 W^2(t)+b_4 W(t)+b_5 G^2(t)+b_5 G(t)W(t)+b_5 G(t)S(t)\\&
+b_6 S(t)+ b_7G(t)+b_8\\
\frac{dS}{dt}&= c_{1}+c_{2}G(t)+c_{3}S(t)+c_{4}W(t-\tau)+c_{5}W^2(t-\tau)+
c_{6}W^3(t-\tau)+c_{7}W^4(t-\tau)\\&
+c_{8}W(t)+c_{9}W^2(t)+c_{10}W^3(t)+c_{11}W^4(t)+c_{12}G^2(t)+h.o.t.,
\end{split}
\end{equation}
The expressions for the coefficients $a_{i},i=1,2,...12; b_{j},j=1,2,...8$; and $c_{k},k=1,2,....12$ are provided in Appendix (\ref{app1}). 
Linearizing (\ref{2.2}) near $(0,0,0)$, equally for the positive equilibrium $E^*=(w^*,g^*,s^*)$ of the original system \ref{2.1}, yields:
\begin{equation}\label{2.4}
\begin{split}
\frac{dW}{dt}&= a_{2}W(t) + a_{3} G(t) + a_{4} S(t-\tau)\\ 
\frac{dG}{dt}&= b_4 W(t)+ b_6 S(t) + b_7 G(t)  \\
\frac{dS}{dt}&=c_{8} W(t) + c_{2} G(t) + c_{3} S(t)+ c_{4} W(t-\tau )
\end{split}
\end{equation}

\section{Existence of equilibria}
To analyze the long-term dynamics of the system (\ref{2.2}), we first establish the existence of a positively invariant and attracting set in the non-negative octant 
$R_{+}^{3}$. This ensures that all trajectories with biologically meaningful initial conditions remain bounded and non-negative for all time. First, the set
\begin{equation*}
\left\lbrace (w,g,s): w(t) \geq 0 , g(t) \geq 0 , s(t) \geq 0 \right\rbrace 
\end{equation*} 
is positively invariant because:
\begin{itemize}
\item When $w=0,\frac{dw}{dt}=0$, $w(t)$ cannot decrease further and remains nonnegative.
\item When $g=0$, $\left.\frac{dg}{dt}\right|_{g=0}=\frac{bw}{1+w} > 0 \quad  if \quad w > 0$, this ensures that $g(t)$ remains nonnegative.
\item When $s=0$, $\left.\frac{ds}{dt}\right|_{s=0}=\frac{cw}{1+w} > 0 \quad  if \quad w > 0$, this ensures that $s(t)$ stays nonnegative.
\end{itemize}
Next, we derive upper bounds for $w(t), g(t)$ and $s(t)$. From the first equation of system (\ref{2.2}), assuming $g=0$ and $s(t-\tau)=0$ we have:
\[ \frac{a w^{2}}{1+w}=\xi_{1}w^{2}  \]
and multiply by $(1+w)$, we get 
\[ a w^{2}= \xi_{1}w^{2}+\xi_{1}w^{3} \]
For large $w$, the dominant term is $\xi_1w^{3}$.Thus
\[ \frac{dw}{dt} \leq  a - \xi_1 w  \]
which implies the upper bound for $w(t)$ is
\[ \limsup_{t\to\infty} w(t) \leq \frac{a}{\xi_1}  \]
From the second equation of system (\ref{2.2}) we have: 
\[\frac{dg}{dt} \leq  b-\xi_2g^{2} \]
Since $\frac{bw}{1+w} \leq b $ and $\xi_{2}(w+g+s)g \geq \xi_{2}^{2} g $,
we have:
\[ \lim_{t\to\infty} \sup g(t) \leq \sqrt{\frac{b}{\xi_2}}  \]
From the third equation of system (\ref{2.2}), assuming $g=0$ and $w(t-\tau)=0$ we have:
\[ \frac{ds}{dt} \leq c + (r -\xi_{3}s) s  \] 
Since $\frac{cw}{1+w} \leq c$ and $\frac{rs+aw(t-\tau)}{1+w(t-\tau)+g+s} \leq r$, for large $s$
we have:
\[ \limsup_{t\to\infty} s(t) \leq \frac{r+\sqrt{r^{2}+4\xi_3 c} }{2\xi_3}  \]
Now, we define the set
\[\Omega= \left\lbrace (w,g,s): 0 \leq w \leq \frac{a}{\xi_1}, 0 \leq g \leq \sqrt{\frac{b}{\xi_{2}}} , 0 \leq s \leq \frac{r+\sqrt{r^{2}+4\xi_3 c} }{2\xi_3} 
\right\rbrace \]
This set $\Omega$ is positively invariant for system (\ref{2.1}), and we only consider $(w,g,s)\in \Omega$. Furthermore, $\Omega$ is attracting as all system trajectories eventually enter and remain within $\Omega$.

Since the delay $\tau$ does not change the number of equilibrium solutions in the system (\ref{2.2}), we analyze the equilibria of the corresponding non-delayed ODE system.
\begin{equation}\label{3.1}
\begin{split}
\frac{dw}{dt}&=\Big(\frac{aw(t)+rs(t)}{1+w(t)+g(t)+s(t)}-(\xi_{1}(w(t)+g(t)+s(t))\Big) w(t)\\ 
\frac{dg}{dt}&=\frac{bw(t)}{1+w(t)}-\Big(\xi_{2}(w(t)+g(t)+s(t))\Big) g(t) \\
\frac{ds}{dt}&=\frac{cw(t)}{1+w(t)}+[\frac{rs(t)+aw(t)}{1+w(t)+g(t)+s(t)}-(\xi_{3}(w(t)+g(t)+s(t))]s(t) 
\end{split}
\end{equation}
The origin $\left(0,0,0\right)$ is the trivial equilibrium of the system (\ref{3.1}). Moreover, we see the boundary equilibrium $E^{0}=(0,0,s^{0})$, where
\begin{equation}\label{3.2}
s^{0}=\frac{r-\xi _3}{\xi _3} ,\quad  r > \xi _3
\end{equation}
This equilibrium represents a scenario where wild and sterile mosquitoes are absent due to high mortality rates $(\xi_1, \xi_2 \gg1 )$. Non-sterile mosquitoes $s$ persist at a level determined by the balance between their intrinsic growth rate $r$ and mortality rate $\xi_{3}$.

Let $N=w+g+s$ denote the total population at equilibrium. A positive equilibrium satisfies 
\begin{equation}\label{3.3}
\left\{
\begin{aligned}
w^{*}&= \frac{
\begin{array}{l}
\left(\xi_1 - \xi_3\right) N \left(\xi_1 N (N+1) - a\right) + c r - \sqrt{\alpha} \\

\end{array}
}{2 a \left(\xi_1 - \xi_3\right) N} \\[10pt]
g^{*} &= \frac{
\begin{array}{l}
b \big( c r - \left(\xi_1 - \xi_3\right) N \left(a + \xi_1 N (N+1)\right)\big)- b \sqrt{\alpha} \\
 
\end{array}
}{2 c \xi_2 N r} \\
s^{*} &= \frac{
\begin{array}{l}
\left(\xi_1 - \xi_3\right) N \left(a + \xi_1 N (N+1)\right) - c r + \sqrt{\alpha}\\

\end{array}
}{2 \left(\xi_1 - \xi_3\right) N r}
\end{aligned}
\right.
\end{equation}

where
$$\alpha=\big(\left(\xi_3 - \xi_1\right) N \left(a + \xi_1 N (N+1)\right) + c r\big)^2 
+ 4 c \xi_1 \left(\xi_1 - \xi_3\right) (N+1) N^2 r$$
$N$ is a positive solution of $H(N)=0$, such that $H(N)=N-(w+g+s)$ of (\ref{3.3}). The positive equilibrium $E^{*}=(w^{*},g^{*},s^{*})$ exists if $w^{*} >0$, $g^{*}>0$ and $s^{*} >0$. Under the assumption that all parameters are positive, this reduces to
\begin{equation*} \label{h1}
\xi _3 > \xi _1    \tag{H1}
\end{equation*}

\section{Stability of the equilibria}
To better understand the system's behavior near the equilibrium points, we apply a Taylor expansion up to the fourth order of the system (\ref{2.2}). The expanded system takes the following form
\begin{equation}\label{4.1}
\begin{split}
\frac{dw}{dt}&=a w^2-a g w^2+a g^2 w^2-a g^3 w^2+a g^4 w^2-a w^3+2 a g w^3-3 a g^2 w^3+4a g^3 w^3 \\&
-5 a g^4 w^3+a w^4-3 a g w^4+6 a g^2 w^4-10 a g^3 w^4+15 a g^4 w^4+r w s(t-\tau )\\&
-g r w s(t-\tau )+g^2 r w s(t-\tau )-g^3 r w s(t-\tau )+g^4 r w s(t-\tau )
-a w^2 s(t-\tau ) \\& 
+2 a g w^2 s(t-\tau )-3 a g^2 w^2 s(t-\tau )+4 a g^3 w^2 s(t-\tau )-5 a g^4 w^2 
s(t-\tau ) \\&
-r w^2 s(t-\tau )+2 g r w^2 s(t-\tau )-3 g^2 r w^2 s(t-\tau )+4 g^3r w^2s(t-\tau )\\&
-5 g^4 r w^2 s(t-\tau )+2 a w^3 s(t-\tau )-6 a g w^3 s(t-\tau )+12 a g^2 w^3
s(t-\tau )+h.o.t., \\
\frac{dg}{dt}&=b w-b w^2+b w^3-b w^4-g^2 \xi _2-g s \xi _2-g w \xi _2   \\
\frac{ds}{dt}&=r s^2-g r s^2+g^2 r s^2-g^3 r s^2+g^4 r s^2-r s^3+2 g r s^3-3 g^2 r s^3+4 g^3 r s^3 \\&
-5 g^4 r s^3+r s^4-3 g r s^4+6 g^2 r s^4-10 g^3 r s^4c w+15 g^4 r s^4-c w^2+c w^3-c w^4 \\&
-g \xi _3 s -s^2 \xi _3 + a s w(t-\tau )-a g s w(t-\tau ) +a g^2 s w(t-\tau )-a g^3 s w(t-\tau ) \\&
+a g^4 s w(t-\tau )-a s^2 w(t-\tau )+2 a g s^2 w(t-\tau )-3 a g^2 s^2 w(t-\tau )-r s^2 w(t-\tau ) \\&
+4 a g^3 s^2 w(t-\tau )-5 a g^4 s^2 w(t-\tau )+2 g r s^2 w(t-\tau )-3 g^2 r s^2 
w(t-\tau ) \\&
+4 g^3 r s^2 w(t-\tau )-5 g^4 r s^2 w(t-\tau )+a s^3 w(t-\tau )-3 a g s^3 w(t-\tau )+h.o.t.,
\end{split}
\end{equation}
To determine the stability of the equilibria, we calculate the characteristic equation 
\begin{equation} \label{4.2}
p\left(\lambda\right):=det\left(\lambda I-J_{0}-J_{1}e^{-\lambda \tau}\right)=0
\end{equation} 
at the equilibria including $(w,g,s)=(0,0,0),(0,0,\frac{r-\xi _3}{\xi _3})$, and $E^*=(w^{*},g^{*},s^{*})$.
If all roots of (\ref{4.2}) have a negative real part, as a direct result of the Grobman-Hartman theorem, the equilibrium solution is locally asymptotically stable \cite{smith2011}.
At the equilibrium $(0,0,0)$, $J_{1}$ is zero matrix, and the Jacobian matrix $J_{0}$ of system (\ref{4.1}) has the form  
\begin{equation*}\label{4.3}
J_{0}=
\begin{pmatrix}
0&0&0\\
b&0&0\\
c&0&0
\end{pmatrix}
\end{equation*}
It follows that the characteristic equation (\ref{4.2}) becomes $p\left(\lambda\right):=\lambda^{3}=0$, has three zero roots 
$\lambda_{1}=\lambda_{2}=\lambda_{3}=0$. So $ E_{0} $ is the center, and the system has no exponential growth or decay modes associated with nonzero eigenvalues. This implies the system states may remain stationary without being attracted or repelled exponentially.

At $ E=(0,0,\frac{r-\xi _3}{\xi _3})$, the Jacobian matrix's $J_{0}$ and $J_{1}$ of system (\ref{4.1}) has the form

\begin{equation*}\label{eq:matrices}
J_{0}=\begin{pmatrix}
-\frac{\left(r-\xi _3\right) \left(r^4-4 \xi _3 r^3+6 \xi _3^2 r^2+\xi _3^3 \left(\xi _1-4 r\right)\right)}{\xi _3^4} & 0 & 0 \\
b & \xi _2 \left(1-\frac{r}{\xi _3}\right) & 0 \\
c & \frac{\xi _3^5-3 r^5+14 \xi _3 r^4-25 \xi _3^2 r^3+20 \xi _3^3 r^2}{\xi _3^4}-7 r & \frac{\left(r-\xi _3\right){}^2 \left(2 \xi _3^2+4 r^2-7 \xi _3 r\right)}{\xi _3^3}
\end{pmatrix},
\end{equation*}
\begin{equation*}
J_{1}=\begin{pmatrix}
0 & 0 & 0\\
0 & 0 & 0\\
-\frac{\left(r-\xi _3\right) \left(r^3 (a+3 r)+\xi _3 \left(\xi _3 \left(\xi _3 \left(-4 a+\xi _3-6 r\right)+2 r (3 a+7 r)\right)-r^2 (4 a+11 r)\right)\right)}{\xi _3^4} & 0 & 0
\end{pmatrix}
\end{equation*}
The characteristic equation (\ref{4.2}) becomes
\begin{equation}  \label{4.4}
p\left(\lambda\right):=\lambda ^3+\delta _1+\delta _2\lambda + \delta _3\lambda^{2}=0
\end{equation} 
Where
\begin{align*}
\delta_{1}=& -2 \xi _1 \xi _2 \xi _3-\frac{4 \xi _2 r^{10}}{\xi _3^8}+\frac{39 \xi _2 r^9}{\xi _3^7}-\frac{170 \xi _2 r^8}{\xi _3^6}+\frac{436 \xi _2 r^7}{\xi _3^5}-\frac{724 \xi _2 r^6}{\xi _3^4}-\frac{4 \xi _1 \xi _2 r^6}{\xi _3^5} \\&
+\frac{803 \xi _2 r^5}{\xi _3^3}+\frac{23 \xi _1 \xi _2 r^5}{\xi _3^4}-\frac{590 \xi _2 r^4}{\xi _3^2}-\frac{54 \xi _1 \xi _2 r^4}{\xi _3^3}+\frac{274 \xi _2 r^3}{\xi _3}+\frac{66 \xi _1 \xi _2 r^3}{\xi _3^2} \\&
-72 \xi _2 r^2-\frac{44 \xi _1 \xi _2 r^2}{\xi _3}+15 \xi _1 \xi _2 r+8 \xi _2 \xi _3 r   \\ 
\delta_{2}=& \xi _1 \xi _2+2 \xi _1 \xi _3+2 \xi _2 \xi _3-\frac{4 r^9}{\xi _3^7}+\frac{35 r^8}{\xi _3^6}-\frac{135 r^7}{\xi _3^5}+\frac{301 r^6}{\xi _3^4}+\frac{\xi _2 r^6}{\xi _3^5}-\frac{423 r^5}{\xi _3^3}-\frac{4 \xi _1 r^5}{\xi _3^4} \\&
-\frac{10 \xi _2 r^5}{\xi _3^4}+\frac{380 r^4}{\xi _3^2}+\frac{19 \xi _1 r^4}{\xi _3^3}+\frac{34 \xi _2 r^4}{\xi _3^3}-\frac{210 r^3}{\xi _3}-\frac{35 \xi _1 r^3}{\xi _3^2}-\frac{55 \xi _2 r^3}{\xi _3^2}+\frac{31 \xi _1 r^2}{\xi _3} \\&
+\frac{45 \xi _2 r^2}{\xi _3}+\frac{\xi _1 \xi _2 r^2}{\xi _3^2}+64 r^2-13 \xi _1 r-17 \xi _2 r-8 \xi _3 r-\frac{2 \xi _1 \xi _2 r}{\xi _3} \\ 
\delta_{3}=& -\xi _1-\xi _2-2 \xi _3+\frac{r^5}{\xi _3^4}-\frac{9 r^4}{\xi _3^3}+\frac{25 r^3}{\xi _3^2}-\frac{30 r^2}{\xi _3}+\frac{\xi _1 r}{\xi _3}+\frac{\xi _2 r}{\xi _3}+15 r  
\end{align*}
Then equation (\ref{4.4}) has three positive roots 
\begin{equation} 
\begin{split}
\lambda_{1}&= \frac{\xi _2 \left( \xi _3- r\right) }{\xi _3}\\
\lambda_{2}&= \frac{\left(r-\xi _3\right){}^2 \left(2 \xi _3^2+4 r^2-7 \xi _3 r\right)}{\xi _3^3} \\
\lambda_{3}&= -\frac{\left(r-\xi _3\right) \left(\xi _1 \xi _3^3+r^4-4 \xi _3 r^3+6 \xi _3^2 r^2-4 \xi _3^3 r\right)}{\xi _3^4}
\end{split}
\end{equation}
If
 \begin{equation}\label{4.5}
	 0<r<\frac{7 \xi _3}{8}-\frac{1}{8} \sqrt{17} \sqrt{\xi _3^2}\quad and \quad \xi _1>\frac{4 \xi _3 r^3-6 \xi _3^2 r^2-r^4+4 \xi _3^3 r}{\xi _3^3}
\end{equation}
 $\lambda_{1} , \lambda_{2}$ and $\lambda_{3} > 0$. So $E$ is the source point. And if
 
\begin{equation}\label{4.6}
\xi _3<r<\frac{1}{8} \sqrt{17} \sqrt{\xi _3^2}+\frac{7 \xi _3}{8}\quad and \quad \xi _1>\frac{4 \xi _3 r^3-6 \xi _3^2 r^2-r^4+4 \xi _3^3 r}{\xi _3^3}
\end{equation}
$\lambda_{1} , \lambda_{2}$ and $\lambda_{3} < 0$. So $E$ is a sink. If one of these roots $>0$, then $E$ is a saddle point. 

\begin{theorem} (Grobman-Hartman) The boundary equilibrium point is
\begin{itemize}
\item[(i)] unstable (repellor) if condition (\ref{4.5}) hold, for any $\tau>0$. 
\item[(ii)] stable (node) if (\ref{4.6}) hold for any $\tau > 0$.
\item[(iii)] Saddle point with index 1 (repel in one direction and attract in two directions) or 2 (repel in two directions and attract in one direction) if one or two of the roots $\lambda_{i}>0,i=1,2,3$, respectively; both are conductive to chaos.
\end{itemize}
\end{theorem} 
 
Now we investigate the stability of the positive equilibrium $E^{*}$. We consider two special cases: 

\begin{enumerate}
\item[(i)] When $\tau = 0$, the characteristic equation (\ref{4.2}) reduces to a polynomial equation. In this case, the stability of $E^{*}$ can be determined using the Routh-Hurwitz criterion that provides necessary and sufficient conditions for all roots to have negative real parts (i.e., local asymptotic stability).
\item[(ii)] For $\tau > 0$, Lemma (\ref{lemm.eq}) states that the number of roots of (\ref{4.2}) with $Re(\lambda)>0$ can change only if a root crosses the imaginary axis. The transversality condition

\[ Re \Big( \frac{d \lambda}{d \tau} \Big)^{-1} > 0  \]
implies that eigenvalues move from the left-half plane to the right-half plane as $\tau$ increases beyond $\tau_{0}$. Hence, if $E^{*}$ is stable at $\tau=0$, it will lose stability at $\tau = \tau_{0}$ (the critical delay value where roots cross the imaginary axis). If no root crosses the imaginary axis for any $\tau > 0$, the equilibrium exhibits absolute stability (i.e., it remains stable for all $\tau > 0$) \cite{Brauer1987}.
\end{enumerate} 
To avoid long-term expressions and easy calculations, we work with the expanded model \ref{2.2}. The associated characteristic equation of (\ref{2.2}) has the following form 
\begin{equation}\label{4.7}
p(\lambda):=\lambda ^3+\delta _1+\delta _2 \lambda+\delta _3 \lambda ^2+\left( \delta _4 \lambda+\delta_{5}\right) e^{-\lambda \tau }+\left(\delta_{6}\lambda+\delta _7\right)e^{-2 \lambda \tau }=0 
\end{equation}
Where:
\begin{align*}
\delta_1=& a_2 b_6 c_2+a_3 b_4 c_3-a_2 b_7 c_3-a_3 b_6 c_8;\quad 
\delta_2=-a_3 b_4+a_2 b_7+a_2 c_3-b_6 c_2+b_7 c_3 \\ 
\delta_3=& -a_2-b_7-c_3 ; \quad
\delta_4= -a_4 c_8 ;\quad \delta_5=-a_4 b_4 c_2-a_3 b_6 c_4+a_4 b_7 c_8 ; \quad 
\delta_6=-a_4 c_4; \\
\delta_7=& a_4 b_7 c_4
\end{align*}

By multiplying $e^{\lambda\tau}$ on both sides of equation \ref{4.7}, we get 
\begin{equation}\label{4.8}
p(\lambda):=\left(\lambda^3+\delta _1+\delta _2 \lambda+\delta_3 \lambda ^2\right) e^{\lambda\tau}+\left(\delta_4 \lambda+\delta_{5}\right)+\left(\delta _6\lambda+\delta_7\right) e^{-\lambda\tau} =0 
\end{equation}

\begin{lemma}\cite{ruan}\label{lemm.eq}
For the transcendental equation
\begin{equation*}
\begin{aligned}
p\left(\lambda,e^{-\lambda\tau_{1}},....,e^{-\lambda\tau_{m}}\right)&=\lambda^{n}+p_{1}^{(0)}\lambda^{n-1}+.....+p_{n-1}^{(0)}\lambda+p_{n}^{(0)} \\& +
\left[p_{1}^{(1)}\lambda^{n-1}+.....+p_{n-1}^{(1)}\lambda+p_{n}^{(1)}\right]
e^{-\lambda \tau_{1}} \\& 
+ ....+\left[p_{1}^{(m)}\lambda^{n-1}+.....+p_{n-1}^{(m)}\lambda+p_{n}^{(m)}\right]
e^{-\lambda \tau_{m}},
\end{aligned} 
\end{equation*}
where $\tau_{i} \geq 0$ $(i=1,2,...,m)$ and $p_{j}^{(i)}(i=1,2,...,m; j=1,2,....,n)$ are constant, as $\left(\tau_{1},\tau_{2},\tau_{3},....,\tau_{m}\right)$ vary, the sum of the orders of the zeros of $P\left(\lambda,e^{-\lambda \tau_{1}},.....,e^{-\lambda \tau_{m}}\right)$ in the open right half plane can change only if a zero appears on or crosses the imaginary axis.
\end{lemma}

As a result, we consider two cases.
\begin{itemize}
\item[\textbf{case(a)}] For $\tau=0$, equation \ref{4.8}  reduces to 
\end{itemize}

\begin{equation}\label{4.9}
p(\lambda):=\lambda^3+\left(\delta _2+\delta_{4}+\delta _6\right)\lambda+
\delta_3 \lambda^2+\left(\delta _1+\delta _5+\delta _7\right) =0 
\end{equation}
A set of necessary and sufficient conditions for all roots of \ref{4.9} to have a negative real part is given by the well-known Routh-Hurwitz criteria in the following form:

\begin{equation*}\label{h2}
\left(\delta _2+\delta_4+\delta_6\right), 
\delta_3 ,\left(\delta _1+\delta _5+\delta _7\right) > 0  , and  
\left(\delta _2+\delta_4+\delta_6\right) \left(\delta_3\right)-
\left(\delta _1+\delta _5+\delta _7\right) > 0 \tag{H2}
\end{equation*}

\begin{itemize}
\item[\textbf{case(b)}] Suppose now that $\tau > 0$. If $\lambda=0$ is a root of (\ref{4.8}), $\delta_1+\delta_5+\delta_7=0$. However $\delta_1+\delta_5+\delta_7\neq0$ since $\delta_1+\delta_5+\delta_7=-a_4 b_4 c_2+a_2 b_6 c_2+a_3 b_4 c_3-a_2 b_7 c_3-a_3 b_6 c_4+a_4 b_7 c_4-a_3 b_6 c_8+a_4 b_7 c_8$ and from (\ref{h2}) $\delta_1+\delta_5+\delta_7 > 0$. Thus $\lambda=0$ is not a root of (\ref{4.8}). Meanwhile, we know that $\lambda=i\omega,(\omega>0)$, is a root of (\ref{4.8}) if and only if $\omega$ satisfies 
\end{itemize}
  
\begin{equation*}
\left(-i\omega ^3+\delta_1+i\delta_2\omega-\delta_3\omega^2\right)e^{i \tau \omega}
+i\omega\delta _4+\delta_{5}+\left(i\delta _6\omega + \delta _7\right)  
e^{-i \tau \omega}=0
\end{equation*}

Splitting the real and imaginary parts, we get

\begin{equation}\label{4.10}
\begin{split}
\left(\delta_1-\delta _3 \omega ^2+\delta _7\right)\cos(\tau \omega )+
\left(\omega ^3-\delta _2 \omega + \delta _6 \omega\right) \sin (\tau  \omega )&= 
-\delta _5 \\
\left(\delta_2 \omega-\omega ^3+\delta _6 \omega \right) \cos (\tau  \omega )+\left(\delta_1-\delta_3 \omega^2-\delta _7\right) \sin(\tau\omega)&=-\delta_4 \omega 
\end{split}
\end{equation}

Following \ref{4.10} we have 

\begin{align}
\cos(\tau\omega )&=\frac{\delta _4 \omega ^2 \left(\omega ^2-\delta _2+\delta _6\right)+\delta _5 \left(\delta _3 \omega ^2-\delta _1+\delta _7\right)}{\left(\omega ^3-\delta _2 \omega \right){}^2+\left(\delta _1-\delta _3 \omega ^2\right){}^2-\delta _6^2 \omega ^2-\delta _7^2} \label{4.11} \\
\sin(\tau\omega)&=\frac{\omega  \left(\delta _3 \delta _4 \omega ^2+\delta _5 \left(\delta _2+\delta _6-\omega ^2\right)-\delta _1 \delta _4-\delta _4 \delta _7\right)}{\left(\omega ^3-\delta _2 \omega \right){}^2+\left(\delta _1-\delta _3 \omega ^2\right){}^2-\delta _6^2 \omega ^2-\delta _7^2} \label{4.12}
\end{align}

Squaring both sides of \ref{4.11} and \ref{4.12}, adding them together and equating the result to 1, $(sin^2(x)+cos^2(x)=1)$, we have:

\begin{equation}\label{4.13}
\omega ^{12}+\gamma _1 \omega ^{10}+\gamma _2 \omega ^8+\gamma _3 \omega ^6+\gamma _4 \omega ^4+\gamma _5 \omega ^2+\gamma _6
\end{equation}
Where
\begin{align*}
\gamma_1=& 2 \delta _3^2-4 \delta _2  \\ 
\gamma_2=& \delta _3^4-4 \delta _2 \delta _3^2-4 \delta _1 \delta _3+6 \delta _2^2+\delta _4^2-2 \delta _6^2 \\ 
\gamma_3=& -4 \delta _2^3+2 \delta _3^2 \delta _2^2-2 \delta _4^2 \delta _2+4 \delta _6^2 \delta _2+8 \delta _1 \delta _3 \delta _2-4 \delta _1 \delta _3^3+2 \delta _1^2-\delta _3^2 \delta _4^2-\delta _5^2-2 \delta _3^2 \delta _6^2-2 \delta _7^2 \\&
+4 \delta _3 \delta _4 \delta _5+2 \delta _4^2 \delta _6\\
\gamma_4=& \delta _2^4+\delta _4^2 \delta _2^2-2 \delta _6^2 \delta _2^2-4 \delta _1 \delta _3 \delta _2^2-4 \delta _1^2 \delta _2+2 \delta _5^2 \delta _2+4 \delta _7^2 \delta _2-4 \delta _3 \delta _4 \delta _5 \delta _2-2 \delta _4^2 \delta _6 \delta _2+\delta _6^4 \\& 
+6 \delta _1^2 \delta _3^2+2 \delta _1 \delta _3 \delta _4^2+\delta _3^2 \delta _5^2+\delta _4^2 \delta _6^2+4 \delta _1 \delta _3 \delta _6^2-2 \delta _3^2 \delta _7^2-4 \delta _1 \delta _4 \delta _5+2 \delta _5^2 \delta _6+2 \delta _3 \delta _4^2 \delta _7 \\ 
\gamma_5=& -4 \delta _3 \delta _1^3+2 \delta _2^2 \delta _1^2-\delta _4^2 \delta _1^2-2 \delta _6^2 \delta _1^2-2 \delta _3 \delta _5^2 \delta _1+4 \delta _3 \delta _7^2 \delta _1+4 \delta _2 \delta _4 \delta _5 \delta _1-2 \delta _4^2 \delta _7 \delta _1-\delta _2^2 \delta _5^2 \\&
-\delta _5^2 \delta _6^2-2 \delta _2^2 \delta _7^2-\delta _4^2 \delta _7^2+2 \delta _6^2 \delta _7^2-2 \delta _2 \delta _5^2 \delta _6+2 \delta _3 \delta _5^2 \delta _7+4 \delta _4 \delta _5 \delta _6 \delta _7 \\
\gamma_6=& \delta _1^4+\delta _5^2 \delta _1^2-2 \delta _7^2 \delta _1^2-2 \delta _5^2 \delta _7 \delta _1+\delta _7^4+\delta _5^2 \delta _7^2
\end{align*}

Let $\omega^2=m$, then equation \ref{4.13} become 

\begin{equation}\label{4.14}
m^{6}+\gamma _1 m ^{5}+\gamma _2 m ^4+\gamma _3 m ^3+\gamma _4 m ^2+\gamma _5 m+\gamma _6 = 0
\end{equation}

Denote 
\begin{equation}\label{4.15}
F(m)= m^{6}+\gamma_1 m ^{5}+\gamma_2 m ^4+\gamma_3 m ^3+\gamma_4 m ^2+\gamma_5 m+\gamma_6 = 0 
\end{equation}

We present the following three lemmas without explicit proof.

\begin{lemma}
If $\gamma_6 < 0$, equation \ref{4.15} has at least one positive root.
\end{lemma}

\begin{lemma}
If $\gamma_6 > 0 $, all roots of the equation \ref{4.15} have negative genuine parts for all delays, and equilibrium $E^{*}$ is locally asymptotically stable.
\end{lemma}

\begin{lemma}
If $\gamma_6 > 0 $ , the necessary condition for equation \ref{4.15} to have positive real root is 
\begin{equation*}\label{h3}
 \delta _5^{2}  < 2\delta _7^{2} \quad and  \quad \left(\delta _1^{4}+\delta _7^{4}+\delta _5^{2}\delta _7^{2}\right) < 2\delta _5^{2}\delta _7\delta _1 \tag{H3}
\end{equation*}
\end{lemma}

Since $F(0) < 0 $, if \ref{h3} holds and $\lim_{m \to +\infty}F(m)=+\infty$, we obtain that \ref{4.14} has at least one positive root and the equation \ref{4.13} will have a positive root. If the equation $F(m)=0$ has more than one positive root, it must have exactly six positive roots. We assume that equation \ref{4.14} has six positive roots, denoted by $m_{1}, m_{2}, m_{3}, m_{4}, m_{5}$ and $m_{6}$, respectively. Then \ref{4.13} has six positive roots 

$\omega_{1}=\sqrt{m_{1}}, \quad  \omega_{2}=\sqrt{m_{2}}, \quad  \omega_{3}=\sqrt{m_{3}}, \quad   \omega_{4}=\sqrt{m_{4}}, \quad   \omega_{5}=\sqrt{m_{5}}, \quad \omega_{6}=\sqrt{m_{6}}$ 
    
By equation (\ref{4.11}), we have 
\begin{equation}
\begin{split}
\tau_k^{(j)}&=\frac{1}{\omega_k}\left[\cos ^{(-1)}\left(\frac{\delta _4 \omega _k^2 \left(\omega _k^2-\delta _2+\delta _6\right)+\delta _5\left( \omega _k^2 \delta _3-\delta _1
+\delta_{7}\right)}{\left(\omega _k^3-\omega _k \delta _2\right){^2}+\left(\delta_{1}-
\omega_k^2 \delta _3 \right){^2}-\delta _6^2 \omega_{k}^2-\delta _7^2}\right)+2 \pi j \right] 
\end{split}
\end{equation}
This implies that the characteristic equation (\ref{4.8}) will have a pair of purely imaginary roots, $\pm i\omega$ with $\tau =\tau_k^{(j)}$, k=1,2,....,6 and j=0,1,2,.....

\[ \lim_{j\to +\infty} \tau_k^{(j)} = +\infty ,k=1,2,...6 \]
Thus, we can define 
\[ \tau_{0}=\tau_{k_{0}}^{(0)} = \min_{1 \leq k \leq 6, j\geq 1} \left\{\tau_k^{(0)}\right\}  , \omega_{0}= \omega_{k0}   \]      
Let $ \lambda(\tau) = \alpha(\tau)+i\omega(\tau)$ be a root of (\ref{4.7}) ,
$\tau=\tau_k^{(j)} , \alpha(\tau_k^{(j)})= 0 $ and  
$\omega(\tau_k^{(j)})= \omega_{k}$

\section{Analysis of Hopf Bifurcation} 
To study the Hopf bifurcation analysis of the system (\ref{2.2}) arising when there is a pair of purely complex roots, and for which we need to verify the transversality condition \cite{hassard1981} as:
\[ \left\{\frac{d ( Re \lambda)}{d\tau} \right\}_{\tau= \tau_{k}^{j}}  > 0  \]
 
This designates at least one eigenvalue with a positive fundamental part for 
$\tau > \tau_{k}$ and preserving the conditions for the existence of a periodic solution. We are interested in finding the purely imaginary roots of the form $\lambda = \pm i\omega$ of the characteristic Eq. (\ref{4.8}). The main goal of this study is to investigate the direction of motion of $\lambda$ when $\tau$ is varied.
Suppose $\lambda_{k}=\alpha_{k}(\tau)+i\omega_{k}(\tau)$ denotes the root of (\ref{4.3}) near $\tau=\tau_{k}$ satisfying $\alpha_{k}(\tau_{k})=0$
and $\omega_{k}(\tau_{k})=\omega_{0},k=1,2,3,....$. Then we have the following theorem;

\begin{theorem}\label{thm:bigthm}
For the characteristic equation (\ref{4.7}), we have the following transversal conditions 

\[ Re\left\{\left(\frac{d\lambda}{d\tau}\right)^{-1}\right\}_{\tau=\tau_k^{(j)}} >0  \]

\begin{proof}
Proof: Substituting $\lambda(\tau)$ into the left-hand side of
(\ref{4.7}) and taking the derivative with respect to $\tau$, we obtain

\begin{align*} 
&\left[3 \lambda^2+\delta_2 + 2 \delta _3 \lambda+\left(\delta_4-\tau \left(  \left(\delta _4 \lambda +\delta _5 \right) e^{-\lambda \tau }
+ \left( 2 \delta _6 \lambda+ 2\delta _7\right)\right) e^{-2 \lambda \tau }\right)\right]\frac{d\lambda}{d\tau}
\\& -\lambda  \left(\left(\delta _4 \lambda + \delta _5\right) e^{-\lambda \tau }+
\left( 2 \delta _6 \lambda + 2 \delta _7\right) e^{-2 \lambda  \tau }\right)=0
\end{align*} 
Consequently,

\begin{align*} 
\left(\frac{d\lambda}{d\tau}\right)^{-1}=&\frac{3\lambda^2+\delta_2+2\delta_3 \lambda+\delta_4 e^{-\lambda\tau} +\delta _6 e^{-2\lambda\tau}}
{\lambda \left[ \left(\delta _4 \lambda +\delta _5 \right) e^{-\lambda \tau }+
\left( 2 \delta _6 \lambda + 2\delta _7\right) e^{-2 \lambda \tau }\right] }-\frac{\tau}{\lambda}
\end{align*} 

Then 
\begin{align*} 
sign\left\{Re\left(\frac{d\lambda}{d\tau}\right)^{-1}\right\}_{\lambda=i\omega_{k}}=\quad &sign\left\{\frac{\left(3\lambda^2+\delta_2+2\delta_3\lambda\right) e^{\lambda\tau} +\delta_4 +\delta_6 e^{-\lambda\tau}}
{\lambda \left(\delta _4 \lambda +\delta _5 +\left(2\delta _6 \lambda + 2 \delta_7\right) e^{-\lambda \tau }\right)} \right\}_{\lambda=i\omega_{k}} \\
=\quad &sign Re\left\{\frac{\left(-3 \omega_{k}^2+\delta_2+2 i \delta _3 \omega_{k}\right) 
e^{i \tau_{k}^{(j)} \omega_{k} }+ \delta _4 +\delta _6 e^{-i \tau_{k}^{(j)}  \omega_{k} }}
{-\delta _4 \omega_{k}^2 + i \delta _5 \omega_{k} - 2 \delta_6 \omega_{k} 
e^{-i \tau_{k}^{(j)} \omega_{k}}  + 2 i \delta _7 \omega_{k}  e^{-i \tau_{k}^{(j)}  \omega_{k} }}\right\}\\
=\quad &sign\left\{\frac{\rho_{1}\rho_{3}+\rho_{2}\rho_{4}}
{\rho_{1}^{2}+\rho_{2}^{2}}\right\}
\end{align*} 

where 
\begin{align*}
\rho_1=& 2\delta_7\omega_{k}\sin(\tau_k^{(j)}\omega_{k})-2 \delta _6 \omega_{k}^2 
\cos(\tau_k^{(j)}\omega_{k})- \delta_{4}\omega_{k}  \\
\rho_2=& 2\delta _6 \omega_{k}\sin(\tau_k^{(j)} \omega_{k} )+2\delta_7 \omega_{k}\cos(\tau_k^{(j)} \omega_{k})+\delta _7 \omega_{k} \\
\rho_3=& \left(\delta_2+3 \delta_{6}-3\omega_{k}^{2}\right)\cos(\tau_k^{(j)}\omega_{k})-2\delta_3\omega_{k} \sin(\tau_k^{(j)}\omega_{k})+\delta _4 \\
\rho_4= & \left(\delta_2-\delta_{6}-3\omega_{k}^{2}\right) \sin(\tau_k^{(j)}\omega_{k})+2\delta_{3}\omega_{k}cos(\tau_k^{(j)}\omega_{k})
\end{align*}
If $\omega \neq 0$, this derivative is positive, so the root must pass from the
negative to the positive half-plane as $\tau$ increases. On the other hand, 
$\omega = 0$ corresponds to a zero root $\lambda=0$, which is impossible since we have assumed that $\delta_{1}+\delta_{4}+\delta_{6} > 0$. Thus, roots can cross the imaginary axis only from left to right as $\tau$ increases. Stability switching may occur, and if it is lost at some critical value of $\tau$ (or does not exist for $\tau = 0$), it can never be regained.

It is necessary to make the following assumption:
\begin{align*}\label{h4} 
\rho_{1}\rho_{3}+\rho_{2}\rho_{4} > 0 \tag{H4}
\end{align*} 
It  means 
\begin{align*}\label{h5} 
Re \left\{\frac{\rho_{3}+i\rho_{4}}{\rho_{1}+i\rho_{3}}\right\} =\left\{\frac{\rho_{1}\rho_{3}+\rho_{2}\rho_{4}}{\rho_{1}^{2}+\rho_{2}^{2}}\right\} > 0 \tag{H5}
\end{align*} 
Therefore, the transversality conditions are satisfied as
\begin{equation}
\eta(\omega_k)= sign\left\{Re\left(\frac{d\lambda}{d\tau}\right)^{-1}\bigm|_{\lambda=i\omega_{k} , \tau=\tau_k^{(j)}} \right\} > 0 
\end{equation}
This completes the proof.
\end{proof}
\end{theorem} 

\begin{theorem}\label{thm:mvt}
Assume that (\ref{h1}) to (\ref{h5}) hold , we get the following result

\begin{itemize}
\item[(i)] The positive equilibrium solution $E^{*}(w^{*},g^{*},s^{*})$ of system (\ref{2.2}) is asymptotically stable for $\tau \in \left[0,\tau_{0}\right]  $ 
\item[(ii)] when  $\tau=\tau_k^{(j)}(k=1,2,...6; j=0,1,2....)$, system (\ref{2.2})  undergoes a Hopf bifurcation at the positive equilibrium $E^{*}(w^{*},g^{*},s^{*})$.
\end{itemize}
\end{theorem}

In theorem \ref{thm:mvt} (ii), the conditions $(H1)$ to $(H5)$ guarantee the occurrence of the Hopf bifurcation at some $\tau_{0}$. However, these conditions are not verified analytically because explicit formulas of $(w^*,g^*,s^*)$ for the general parameters in \ref{2.1} are long expressions. Below, we show that the Hopf bifurcations do occur for some parameter ranges.

\newpage

\section{Numerical simulation} 
In this section, we perform numerical simulations to validate the analytical results and investigate the impact of time delay $\tau$ on the system dynamics (\ref{2.2}). The parameter values in Table (\ref{tab}) are chosen within biologically plausible ranges derived from empirical studies.

\begin{table}[H]\label{tab}
\footnotesize
\begin{center}
  \begin{tabular}{|c|c|c|c|} \hline
   Par. & \bf Definition & \bf Range & \bf Unit  \\ \hline
    $a$ & \parbox{4cm}{wild offspring production rate} &[1 ,20]& \parbox{4cm}{offspring mosquito per day} \\  \hline
    $b$ & \parbox{5cm}{release rate of the sterile mosquitoes} &[5,25]& \parbox{4cm}{sterile mosquitoes released per day} \\ \hline
     $c$ & \parbox{5cm}{residual fertility or release rate of the non-sterile mosquitoes} &[0.01, 25]& 
    {non-sterile mosquitoes released per day} \\ \hline
    $r$ & \parbox{5cm}{non-sterile insects offspring production rate } &[1,20]& {offspring mosquito per day} \\ \hline
    $\xi_{1}$& \parbox{5cm}{density dependent death rates of the wild mosquitoes} &[0.02,3.5] & deaths per mosquito per day \\ \hline
    $\xi_{2}$& \parbox{5cm}{density dependent death rates of the sterile mosquitoes} &[1.5,2.5] &deaths per mosquito per day \\ \hline
    $\xi_{3}$& \parbox{5cm}{density dependent death rates of the non-sterile mosquitoes} & [0.1,2.5] &deaths per mosquito per day \\ \hline
    $\tau$ & delay time & [0.1 ,7] & days\\ \hline
    $w$ &  wild mosquitoes & - &  mosquitoes abundance\\ \hline
    $g$ & sterile mosquitoes & - &  mosquitoes abundance \\ \hline
    $s$ & non-sterile mosquitoes & - & mosquitoes abundance \\ \hline
   \end{tabular}
\end{center}
\caption{Parameters and variables used in numerical analysis of the model (\ref{2.2})}\label{tab1}
\end{table}

\begin{figure}[H]
\centering
\includegraphics[width=0.3\textwidth]{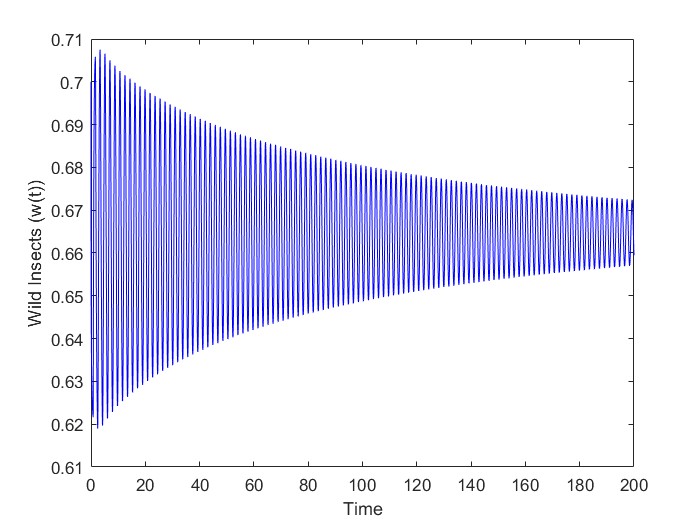}
\includegraphics[width=0.3\textwidth]{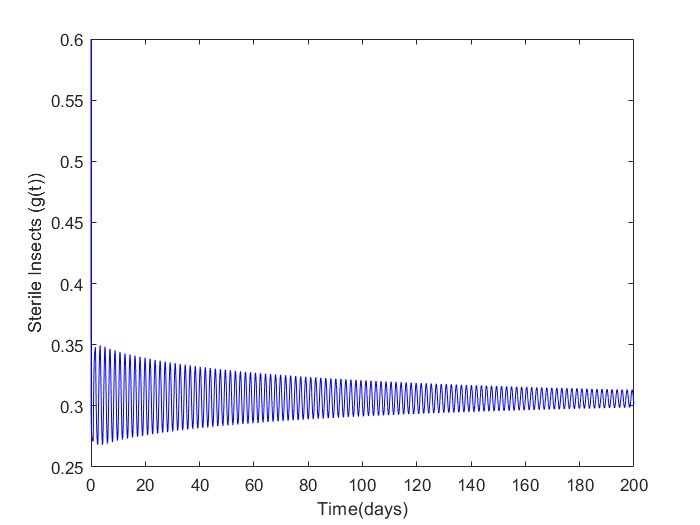}
\includegraphics[width=0.3\textwidth]{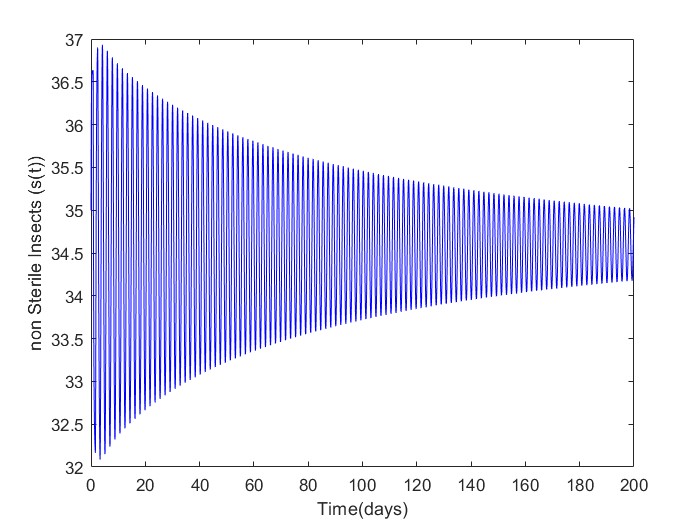}\\
\includegraphics[width=0.3\textwidth]{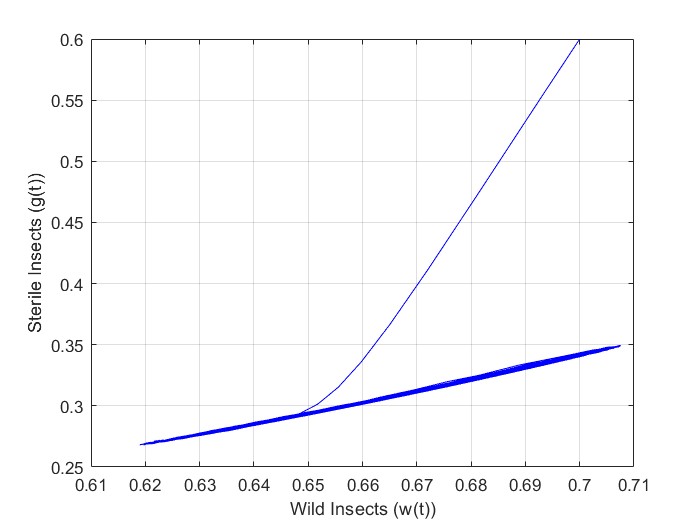}
\includegraphics[width=0.3\textwidth]{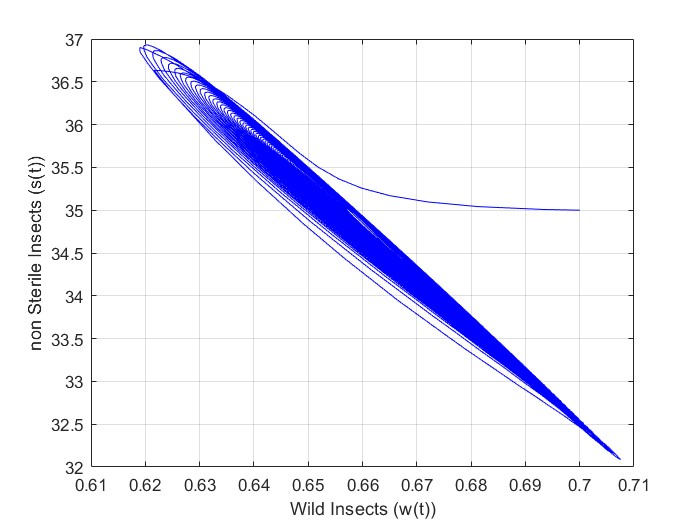}
\includegraphics[width=0.3\textwidth]{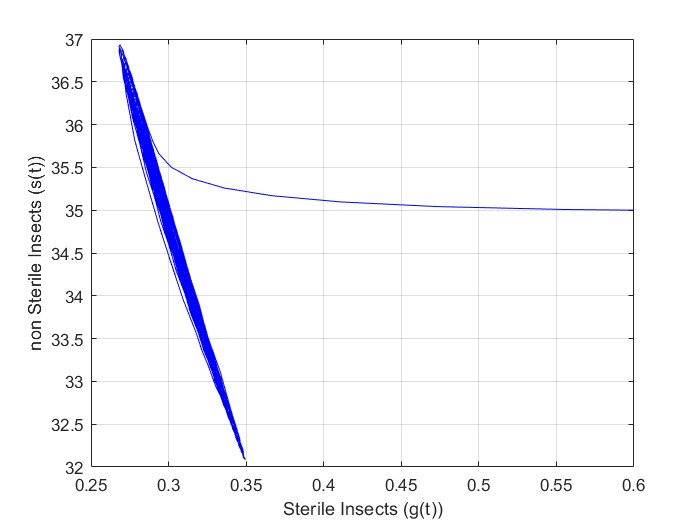}\\
\includegraphics[width=0.3\textwidth]{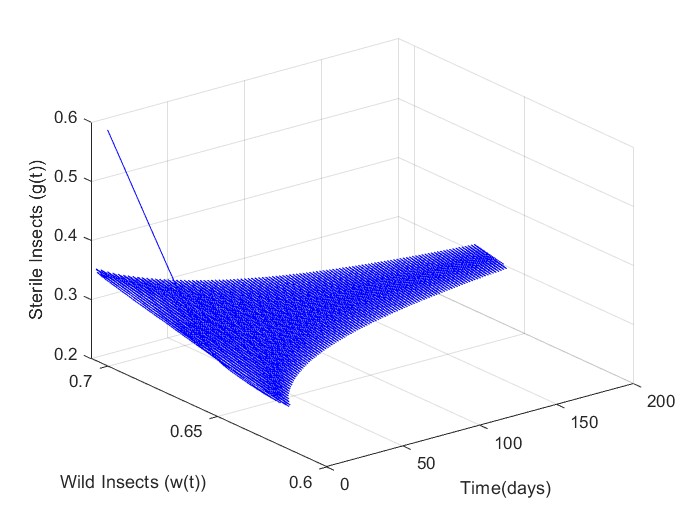}
\includegraphics[width=0.3\textwidth]{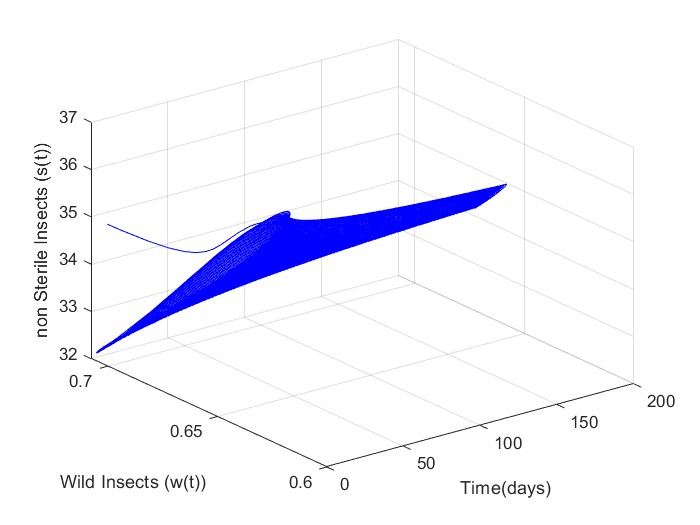}
\includegraphics[width=0.3\textwidth]{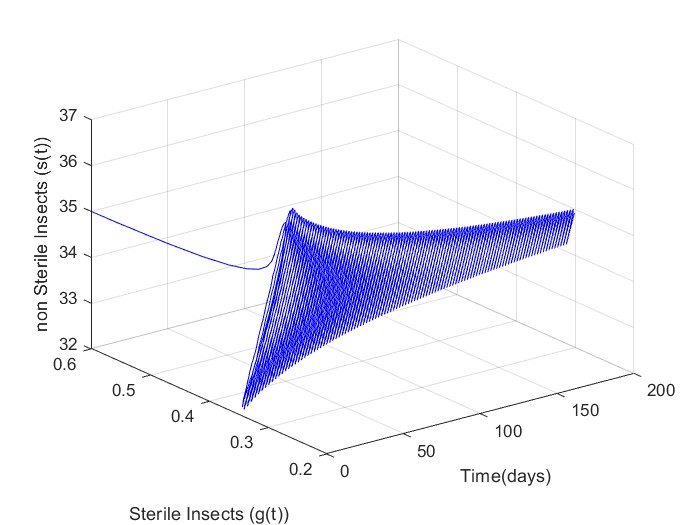}\\
\includegraphics[width=0.3\textwidth]{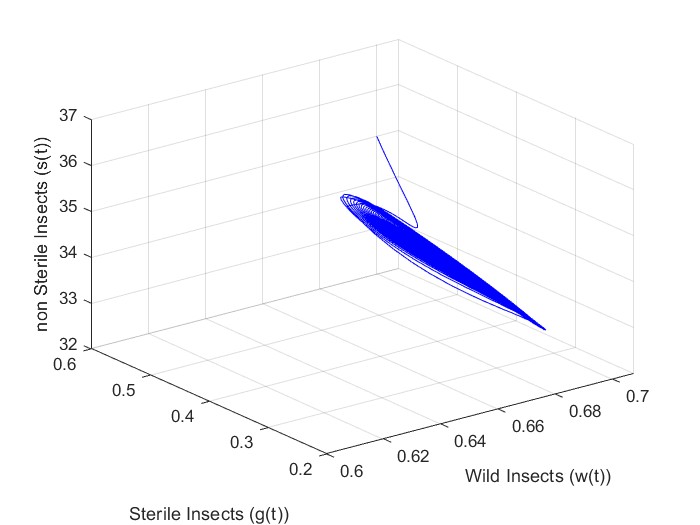}
\caption{Dynamical behaviour and phase portrait of system (\ref{2.2}) with $\tau=0.7 < \tau_{0}\approx  1.74066 $. The positive equilibrium $E^{*}$ is asymptotically stable. The initial value $(w,g,s)=(18.001,0.007,0.005)$.}\label{fig1} 
\end{figure}

\begin{figure}[H]
\centering
\includegraphics[width=0.3\textwidth]{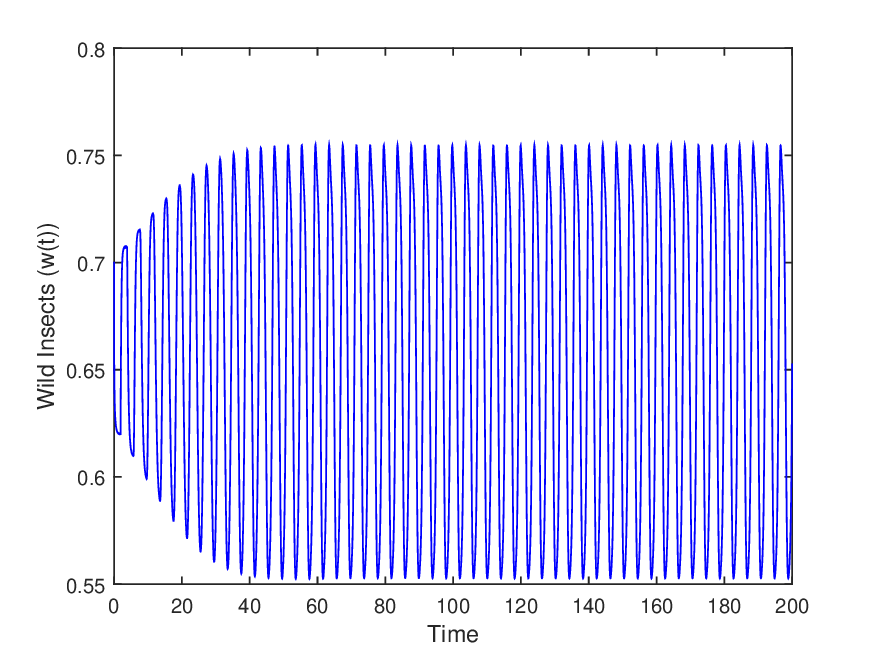}
\includegraphics[width=0.3\textwidth]{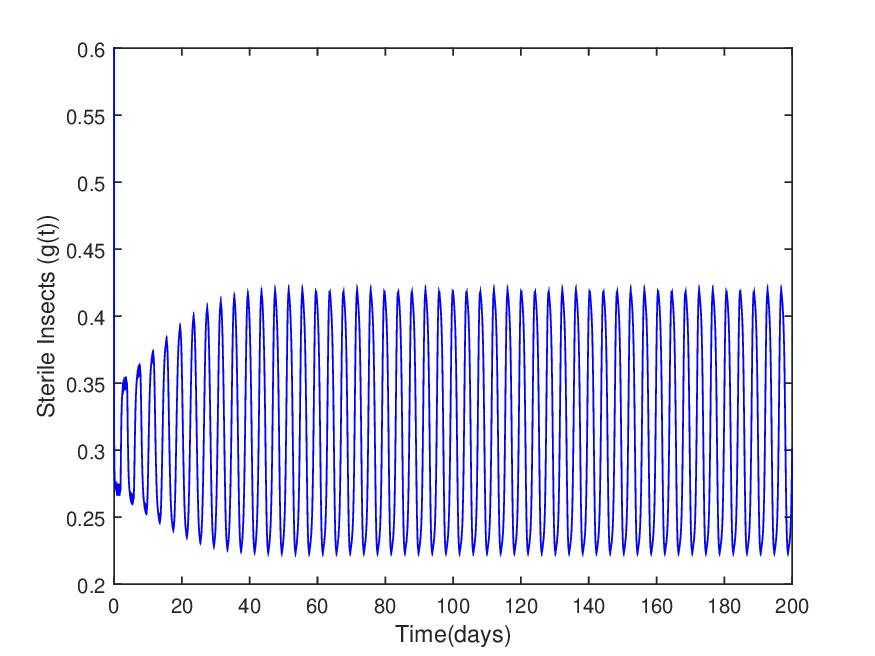}
\includegraphics[width=0.3\textwidth]{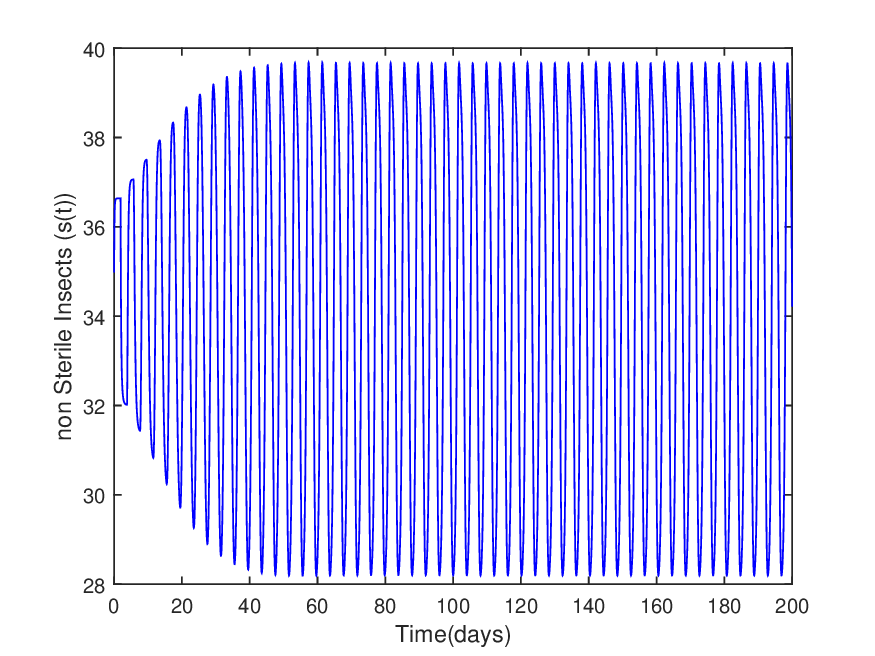}\\
\includegraphics[width=0.3\textwidth]{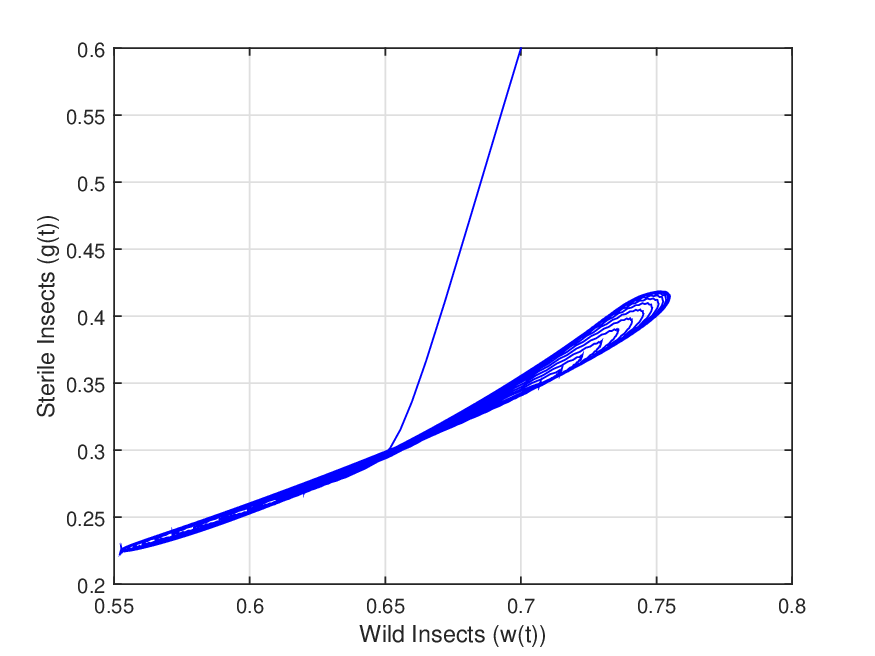}
\includegraphics[width=0.3\textwidth]{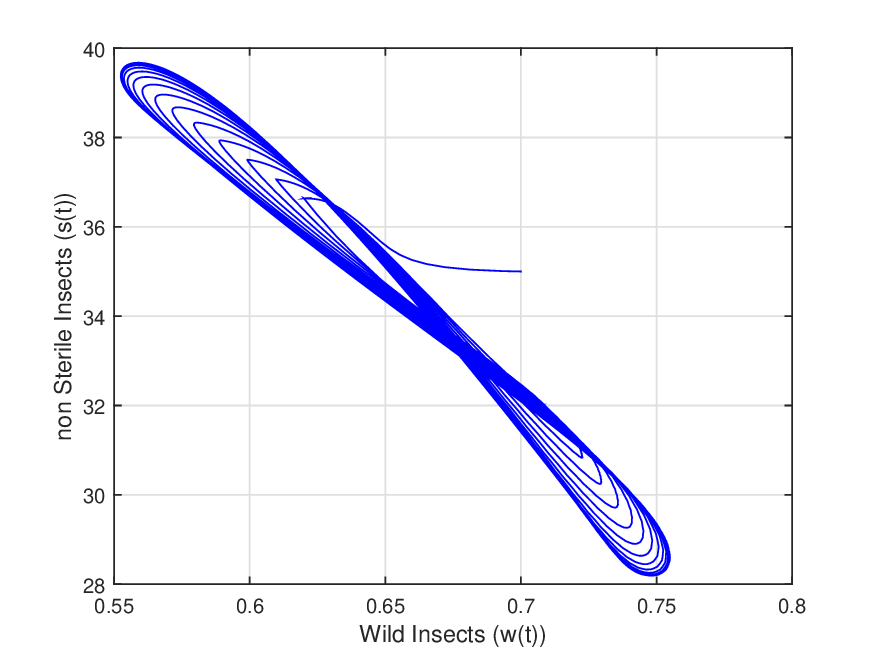}
\includegraphics[width=0.3\textwidth]{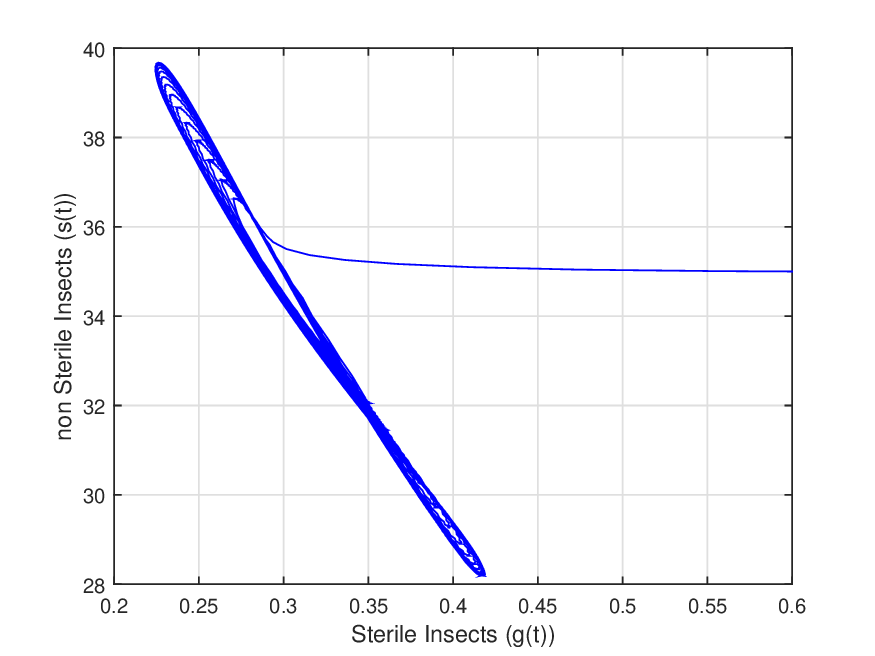}\\
\includegraphics[width=0.3\textwidth]{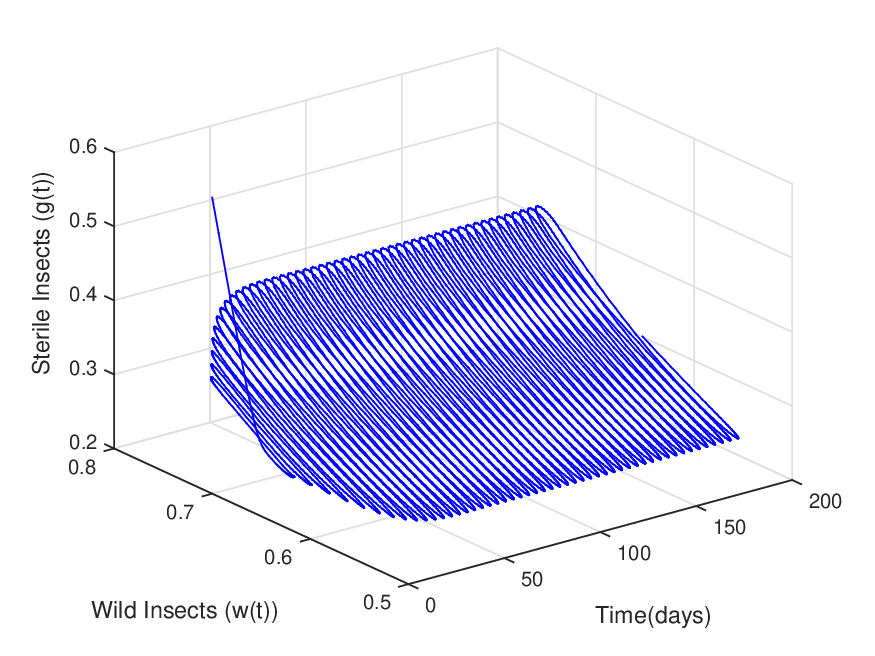}
\includegraphics[width=0.3\textwidth]{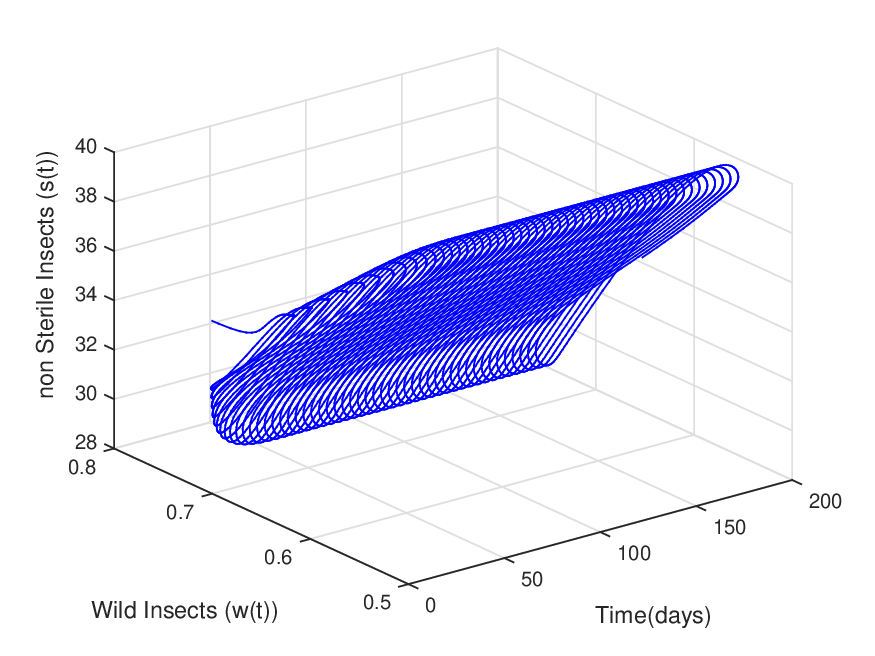}
\includegraphics[width=0.3\textwidth]{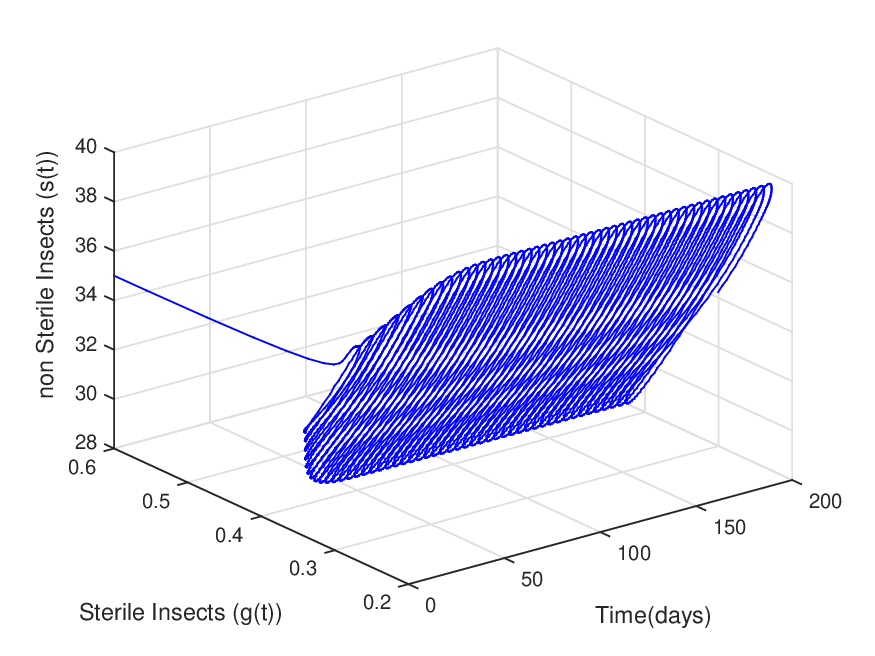}\\
\includegraphics[width=0.3\textwidth]{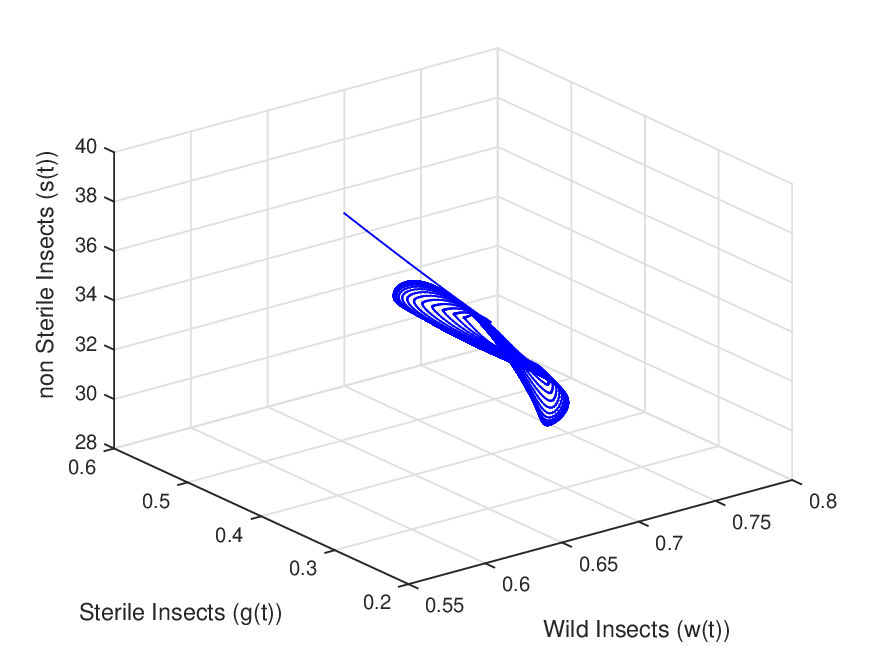}
 \caption{Dynamical behaviour and phase portrait of system (\ref{2.2}) with $\tau =1.9 > \tau_{0} \approx 1.74066$.The origin loses its stability, and Hopf bifurcation occurs. The initial value $(w,g,s)=
 (18.001,0.007,0.005)$.}
 \label{fig2}
\end{figure}

We simulate the system dynamics for varying delay values $\tau$ (see Figure \ref{fig1} and Figure \ref{fig2} ). In Figure \ref{fig1}, stability switch occurs at $\tau=\tau_{0}$. The positive equilibrium $E^{*}$ is stable for $\tau < \tau_{0}$. In Figure \ref{fig2}, when $ \tau $ passes through the critical value $\tau_{0}$, the positive equilibrium $E^{*}$ loses its stability and Hopf bifurcation occurs.

\begin{figure}[H]
\centering
\includegraphics[width=0.96\textwidth]{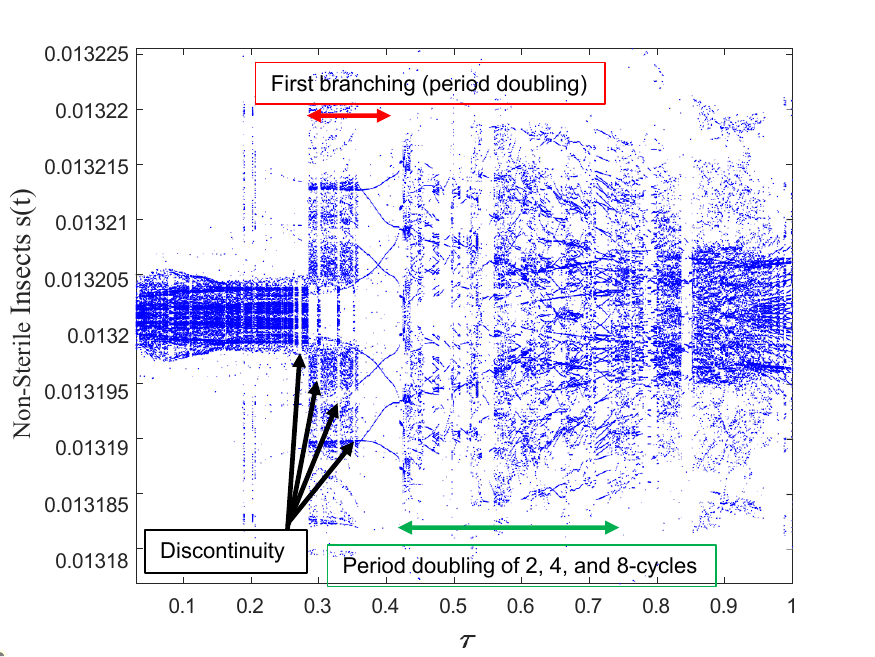}\\
 \caption{One-parameter bifurcation diagram of the non-sterile insects dynamics in system (\ref{2.2}) with	$\tau_{0} \approx0-1$. The initial value is $(w,g,s)=(0.8,0.7,0.6)$ and parameters are taken as $a=5;  b=18;c=0.05; r=1; \xi_1=0.5;\xi_2=0.2;  \xi_3=0.3$.}
\label{fig3} 
\end{figure}

In order to illustrate the change of regimes along the delay values, a one-parameter bifurcation diagram is plotted along the $\tau$ values considering the evolution of non-sterile insect dynamics. Points of variable $s$ are plotted versus $\tau$ , as $\tau$ is increased from $0$ to $1$ (Fig \ref{fig3}). Each vertical slice at a given $\tau$ shows the long-term states of the system after transients die out. For $\tau \approx 0.1-0.25$, there is a stable steady state and the population settles to a nearly constant value. Therefore, the insect population is regulated and predictable. Around $\tau \approx 0.3-0.4$, the first branch splitting (period-doubling bifurcation) occurs. Here, the population stops being constant and begins oscillating in cycles, and wild and non-sterile insects alternately increase/decrease over time. As $\tau$ increases further to around $0.4–0.7$, more period doublings of two-, four-, and eight-cycles occur. Repeated discontinuity occurrence in that region is due to more stable periodic solutions for delays in this region. From this region of delay values, the insect dynamics become increasingly complex with alternating booms and crashes in population size. However, beyond $\tau\approx0.7$, the diagrams become very dense and scattered, indicating chaotic dynamics and subsequently, populations fluctuate in an irregular, unpredictable way. Here, the dynamics are highly sensitive to initial conditions.

\begin{figure}[H]
	\centering
	\includegraphics[width=0.96\textwidth]{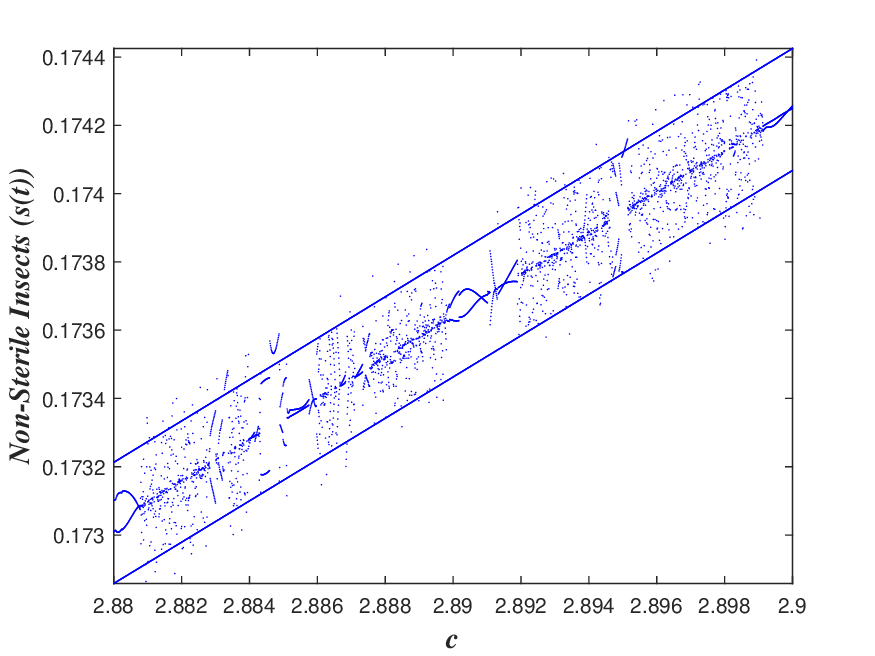}\\
	\caption{One-parameter bifurcation diagram of the non-sterile insects dynamics in system (\ref{2.2}) with	$c\approx2.88-2.90$. The initial value is $(w,g,s)=(0.7,0.8,0.6)$ and parameters are taken as $a=5; b=18; r=1; \xi_1=1.5;  \xi_2=1.2; \xi_3=2.3$.}
	\label{fig4} 
\end{figure}

To reveal the effect of residual fertility on the dynamics of wild and non-sterile insects, a one-dimensional bifurcation diagram was plotted along the c values, seeing non-sterile insects' dynamics (Fig. \ref{fig4}). The diagram shows a clear band of non-sterile insect population values as residual fertility varies within a very narrow range (approximately 2.88 to 2.89). This suggests that within this interval, the system exhibits sensitive dependence on the fertility parameter, likely due to underlying nonlinear dynamics. The scattered, almost fractal-like arrangement of points within the band hints at complex or possibly chaotic dynamics for these values of residual fertility. Rather than converging to a single population level (as seen in simpler, stable regimes), the system exhibits multiple or dense population levels for a given c. Despite the complexity, the non-sterile insect populations remain bounded within a relatively narrow vertical interval (approximately 0.1729 to 0.1738), indicating that the overall population does not experience explosive growth or collapse for the range of residual fertility considered. The diagram suggests that small changes in residual fertility can lead to large qualitative changes in population behavior. This is typical near bifurcation points, where the system may shift from stable, periodic, or simple behavior to more complex or chaotic dynamics.

\begin{example}
Displaying the positive equilibria of \ref{3.1} with the parameters taken from table \ref{tab} 
\begin{equation}
a=18, \quad b=35 ,\quad \xi_1=0.02 ,\quad  \xi_2=1.5, \quad  \xi_3=0.1 , \quad  r=0.99 , \quad  c=0.19 
\end{equation}

Condition (\ref{h1}) is satisfied.The positive solution of $H(N)=0$ is $N=898.972$.The  positive equilibria $E^{*}=(898.943, 0.0259267, 0.00263897)$ is asymptotically stable for system \ref{3.1}.
\end{example}

\section{Discussion}\label{sec12}
In this paper, we developed a three-dimensional model and studied the dynamics of interactive systems involving wild, sterile, and non-sterile mosquitoes with a saturated release rate strategy using delayed ODEs. This was done by incorporating a delay in the growth stage of wild and fertilized mosquito populations in the delayed model (\ref{2.2}). We first established the positivity of solutions, boundedness of the system, and provided a complete stability analysis by considering the case when $\tau=0$ and established stability conditions for the positive equilibrium using the well-known Routh-Hurwitz criterion. For 
$\tau > 0$, we demonstrated that a stability switch occurs via a Hopf bifurcation, where the stable equilibrium becomes unstable. 
We have showed that if the conditions (\ref{h1}),(\ref{h2}) and (\ref{h4}) hold, the positive equilibrium $E(w^{*}(t),g^{*}(t),s^{*}(t))$ of system (\ref{2.2}) is asymptotically stable for all $\tau \in [0, \tau_{0})$. This suggests that the densities of wild, sterile, and non-sterile populations will tend to stabilize. It is believed that delays destabilize the system (e.g., see \cite{cookeVanderish86}). By considering the delay as a bifurcation parameter, we prove the existence of a Hopf bifurcation. Accordingly, we showed that when $\tau$ is small, insect populations remain stable, while a critical delay $\tau_{c}\approx 0.3–0.4$ triggers oscillations and Hopf or period-doubling bifurcation occurs. Longer delays cause unpredictable, chaotic fluctuations in wild, sterile, and non-sterile insects that are not studied analytically here and can be a good subject for future studies. This means that introducing time lags (e.g., in sterile/non-sterile insect release, or non-sterile insect coupling with wilds) can shift the system from stability to periodic oscillations and progressively to chaos.  Additionally, we investigated the local stability of the positive equilibrium $E(w^{*}(t),g^{*}(t),s^{*}(t))$ and analyzed the local Hopf bifurcation of the delayed model.

 Here, for the first time, we have investigated the population with residual fertility (denoted by $s(t)$) by extending a commonly used two-dimensional model through the addition of a new equation.\cite{cai2} showed that a constant release rate of sterile mosquitoes stabilizes the system, whereas a proportional release rate can lead to complex dynamics and oscillatory interactions between wild and sterile populations. As the time delay varies, Hopf and other bifurcations may arise, leading to periodic coexistence between wild and sterile mosquitoes. However, the later workers suggested that oscillations are not caused by simple developmental delays alone, but by more complex delayed processes. These findings highlight that sterile mosquito release strategies play a crucial role in controlling the dynamics \cite{cai2}. \cite{hui} developed a delay differential model for wild and sterile mosquitoes by incorporating maturation delay. For constant sterile releases, the later authors established a release threshold that determines the number and stability of equilibria and proved conditions for global asymptotic stability. Their findings suggest that key dynamical features, such as bistability, are consistent with existing simpler models, indicating that these behaviors may be biologically essential and model-independent. The bifurcation diagram of non-sterile insect dynamics over residual fertility parameter $c$ in the current study visually demonstrates how the non-sterile insect population responds to variation in this parameter. It reveals transitions into complex or chaotic dynamics for certain fertility levels, emphasizing the need for careful parameter management in biological control scenarios and highlighting the rich dynamics possible in ecological models. For vector-borne disease control and pest management, understanding these bifurcations is crucial. Near certain values of residual fertility, population control via fertility reduction might have unpredictable effects, potentially leading to population fluctuations even if average fertility is lowered. The persistence of non-sterile insects within a narrow population range suggests that the ecosystem or model considers density-dependent regulation, resource limitations, or other stabilizing factors. Finally, we suggest studying the zero-Hopf bifurcation of the current model to depict richer dynamics in non-sterile insects, as their stability and success in the field after releasing them are expected; however, their population size and dynamics over time are critical in terms of either disease control or pest management in ecosystems, in case the model is implemented for agricultural insect pests.

\appendix
\section{Expressions for the coefficients in (\ref{2.2})}\label{app1}
In the following expressions, $A_1, A_2, A_3$ are $w^*, g^*$, and $s^*$, respectively.
\begin{align*}
a_1=&\frac{a A_1^2}{1+A_1+A_2+A_3}+\frac{A_3 A_1 r}{1+A_1+A_2+A_3}-A_1^2\xi_1-A_2 A_1 \xi _1-A_3 A_1 \xi_1 \\
a_2=&-\frac{a A_1^2}{\left(1+A_1+A_2+A_3\right)^2}+\frac{2 a A_1}{1+A_1+A_2+A_3}-\frac{A_3 A_1 r}{\left(1+A_1+A_2+A_3\right)^2} \\&
+\frac{A_3 r}{1+A_1+A_2+A_3} 
-2 A_1 \xi _1-A_2 \xi _1-A_3 \xi _1\\
a_3=& -\frac{a A_1^2}{\left(1+A_1+A_2+A_3\right)^{2}}-\frac{A_3 A_1 r}{\left(1+A_1+A_2+A_3\right)^{2}}-A_1 \xi _1\\
a_4=& -\frac{a A_1^2}{\left(1+A_1+A_2+A_3\right)^{2}}+\frac{A_1 r}{1+A_1+A_2+A_3}-\frac{A_3 A_1 r}{\left(1+A_1+A_2+A_3\right)^{2}}-A_1 \xi _1\\
a_5=& \frac{a A_1^2}{\left(1+A_1+A_2+A_3\right)^{3}}-\frac{2 a A_1}{\left(1+A_1+A_2+A_3\right)^{2}}+\frac{a}{1+A_1+A_2+A_3} \\&
+\frac{A_3 A_1 r}{\left(1+A_1+A_2+A_3\right)^{3}}
-\frac{A_3 r}{\left(1+A_1+A_2+A_3\right)^{2}}-\xi _1\\
a_6=&-\frac{a A_1^2}{\left(1+A_1+A_2+A_3\right)^{4}}+\frac{2 a A_1}{\left(1+A_1+A_2+A_3\right)^{3}}-\frac{a}{\left(1+A_1+A_2+A_3\right)^{2}}\\&
-\frac{A_3 A_1 r}{\left(1+A_1+A_2+A_3\right)^{4}}+\frac{A_3 r}{\left(1+A_1+A_2+A_3\right)^{3}} \\
a_7=& \frac{a A_1^2}{\left(1+A_1+A_2+A_3\right)^{5}}-\frac{2 a A_1}{\left(1+A_1+A_2+A_3\right)^{4}}+\frac{a}{\left(1+A_1+A_2+A_3\right)^{3}}\\&
+\frac{A_3 A_1 r}{\left(1+A_1+A_2+A_3\right)^{5}}-\frac{A_3 r}{\left(1+A_1+A_2+A_3\right)^{4}} \\
a_8=&\frac{a A_1^2}{\left(1+A_1+A_2+A_3\right)^{3}}+\frac{A_3 A_1 r}{\left(1+A_1+A_2+A_3\right)^{3}} \\
a_9 =&-\frac{a A_1^2}{\left(1+A_1+A_2+A_3\right)^{4}}-\frac{A_3 A_1 r}{\left(1+A_1+A_2+A_3\right)^{4}}\\
a_{10} =& \frac{a A_1^2}{\left(1+A_1+A_2+A_3\right)^{5}}+\frac{A_3 A_1 r}{\left(1+A_1+A_2+A_3\right)^{5}}\\
a_{11}=& \frac{a A_1^2}{\left(1+A_1+A_2+A_3\right)^{3}}-\frac{A_1 r}{\left(1+A_1+A_2+A_3\right)^{2}}+\frac{A_3 A_1 r}{\left(1+A_1+A_2+A_3\right)^{3}}\\
a_{12}=& -\frac{a A_1^2}{\left(1+A_1+A_2+A_3\right)^{4}}+\frac{A_1 r}{\left(1+A_1+A_2+A_3\right)^{3}}-\frac{A_3 A_1 r}{\left(1+A_1+A_2+A_3\right)^{4}} \\
b_1= & \frac{A_1 b}{\left(1+A_1\right)^{5}}-\frac{b}{\left(1+A_1\right)^{4}} ,\quad  
b_2=\frac{b}{\left(1+A_1\right)^{3}}-\frac{A_1 b}{\left(1+A_1\right)^{4}}, \\
b_3=&\frac{A_1 b}{\left(1+A_1\right)^{3}}-\frac{b}{\left(1+A_1\right)^{2}} \quad
b_4=\frac{b}{1+A_1}-\frac{A_1 b}{\left(1+A_1\right)^{2}}-A_2 \xi _2 \\
\end{align*}
\begin{align*}
b_5=&-\xi _2, \quad
b_6=-A_2 \xi _2 , \quad
b_7=-A_1 \xi _2-2 A_2 \xi _2-A_3 \xi _2, \\
b_8=&\frac{A_1 b}{1+A_1}-A_2^2 \xi _2-A_1 A_2 \xi _2-A_3 A_2 \xi _2 \\ 
c_1= & \frac{a A_1 A_3}{1+A_1+A_2+A_3}+\frac{A_3^2 r}{1+A_1+A_2+A_3}+\frac{A_1 c}{1+A_1}-A_3^2 \xi _3-A_1 A_3 \xi _3-A_2 A_3 \xi _3 \\
c_2=& -\frac{a A_1 A_3}{\left(1+A_1+A_2+A_3\right)^{2}}-\frac{A_3^2 r}{\left(1+A_1+A_2+A_3\right)^{2}}-A_3 \xi _3 \\
c_3= & -\frac{a A_1 A_3}{\left(1+A_1+A_2+A_3\right){^2}}+\frac{a A_1}{1+A_1+A_2+A_3}-\frac{A_3^2 r}{\left(1+A_1+A_2+A_3\right){^2}}
\\& +\frac{2 A_3 r}{1+A_1+A_2+A_3} 
-2 A_3 \xi _3-A_1 \xi _3-A_2 \xi _3 \\
c_4=& \frac{a A_3}{1+A_1+A_2+A_3}-\frac{a A_1 A_3}{\left(1+A_1+A_2+A_3\right){^2}}-A_3 \xi _3-\frac{A_3^2 r}{\left(1+A_1+A_2+A_3\right){^2}}\\
c_5=& -\frac{a A_3}{\left(1+A_1+A_2+A_3\right){^2}}+\frac{a A_1 A_3}{\left(1+A_1+A_2+A_3\right){^3}}+\frac{A_3^2 r}{\left(1+A_1+A_2+A_3\right){^3}}\\
c_6=& \frac{a A_3}{\left(1+A_1+A_2+A_3\right){^3}}-\frac{a A_1 A_3}{\left(1+A_1+A_2+A_3\right){^4}}-\frac{A_3^2 r}{\left(1+A_1+A_2+A_3\right){^4}}\\
c_7=& -\frac{a A_3}{\left(1+A_1+A_2+A_3\right){^4}}+\frac{a A_1 A_3}{\left(1+A_1+A_2+A_3\right){^5}}+\frac{A_3^2 r}{\left(1+A_1+A_2+A_3\right){^5}}\\
\end{align*}
\begin{align*}
c_8=& \frac{c}{1+A_1}-\frac{A_1 c}{\left(1+A_1\right){^2}} , \quad 
c_9= \frac{A_1 c}{\left(1+A_1\right){^3}}-\frac{c}{\left(1+A_1\right){^2}} \\
c_{10}=& \frac{c}{\left(1+A_1\right){^3}}-\frac{A_1 c}{\left(1+A_1\right){^4}} , \quad
c_{11}=\frac{A_1 c}{\left(1+A_1\right){^5}}-\frac{c}{\left(1+A_1\right){^4}} \\
c_{12}=& \frac{a A_1 A_3}{\left(1+A_1+A_2+A_3\right){^3}}+\frac{A_3^2 r}{\left(1+A_1+A_2+A_3\right){^3}}
\end{align*}


\noindent \noindent {\small {\bf Ruqaya Hussein}, Department of Mathematics, Faculty of Mathematical Sciences, Tarbiat Modares
University, Tehran 14115-134, Iran (r.hussein@modares.ac.ir)}
\\
{\small {\bf Sergey Kryzhevich}, Gdańsk University of Technology,ul. Gabriela Narutowicza 11/12, 80-233, Gdańsk, Poland (serkryzh@pg.edu.pl kryzhevicz@gmail.com)}
\\
{\small {\bf Khosro Tajbakhsh}, Department of Mathematics, Faculty of Mathematical Sciences, Tarbiat Modares
University, Tehran 14115-134, Iran (khtajbakhsh@modares.ac.ir)}

 \end{document}